\documentclass[12pt]{article}
\usepackage{latexsym,amsfonts,amssymb}
\usepackage[dvips]{color}
\usepackage{colordvi,multicol}
\usepackage{amsmath}
\usepackage{graphicx}
\usepackage{pdflscape}
\usepackage{rotating}

\allowdisplaybreaks

\usepackage{latexsym,amsfonts,amssymb}
\usepackage[dvips]{color}
\usepackage{colordvi,multicol}
\usepackage{amsmath}
\usepackage{graphicx}
\usepackage{pdflscape}
\usepackage{rotating}
\usepackage{breqn}
\usepackage{mathtools}

\newtheorem{thm}{Theorem}[section]

\newtheorem{pro}[thm]{Proposition}

\newtheorem{rem}[thm]{Remark}

\def\E {{\mathbf E}}
\def\P{{\mathbf P}}
\def\Var{{\rm Var}}

\date{}
\begin{document}

\title{\bf On the anti-concentration functions of some familiar families of distributions}

\author{ Ze-Chun Hu$^{1,}$\thanks{{Corresponding
author.} \vskip 0cm {\it \ E-mail addresses:} zchu@scu.edu.cn (Z.C.
Hu), rsong@illinois.edu (R. Song), 630769091@qq.com (Y. Tan).} , Renming Song$^2$ and Yuan Tan$^{1}$\\ \\
 {\small $^1$College of Mathematics, Sichuan  University,
 Chengdu 610065, China}\\ \\
  {\small $^2$Department of Mathematics,
University of Illinois Urbana-Champaign, Urbana, IL 61801, USA}}

 \date{}
\maketitle

\begin{abstract}
Let $\{X_{\alpha}\}$ be a family of random variables following a certain type of distributions with finite expectation $\E[X_{\alpha}]$
and finite variance $\Var(X_{\alpha})$, where $\alpha$ is a parameter. Motivated by the recent paper of Hollom and Portier (arXiv: 2306.07811v1), we study the anti-concentration function
$(0, \infty)\ni y\to \inf_{\alpha}\P\left(|X_{\alpha}-\E[X_{\alpha}]|\geq y \sqrt{\Var(X_{\alpha})}\right)$ and find its explicit expression.
We show that, for certain familiar families of distributions, including uniform distributions, exponential distributions, non-degenerate Gaussian distributions and student's $t$-distribution, the anti-concentration function is not identically zero, while for some other familiar families of distributions, including binomial, Poisson,
negative binomial, hypergeometric, Gamma, Pareto, Weibull, log-normal and Beta distributions, the anti-concentration function is identically zero.
\end{abstract}

\noindent  {\it MSC:} 60E07,  60E15.

\noindent  {\it Keywords:} Distribution, measure anti-concentration.

\section{Introduction}\setcounter{equation}{0}

Let $\mathcal{R}$ be the family of random variables of the form $X=\sum^n_{k=1}a_k\varepsilon_k$, where $n\ge 1$,  $a_k, k=1, \dots, n,$ are real numbers with $\sum^n_{k=1} a_k^2=1$, and $\varepsilon_k$, $k=1, 2, \dots$, are independent
Rademacher random variables (i.\/e., $\P(\varepsilon_k=1)=\P(\varepsilon_k=-1)=1/2$).

Tomaszewski's conjecture, see \cite{Guy}, can be stated as $\P(|X|\leq 1)\geq 1/2$ for all $X\in\mathcal{R}$. Tomaszewski's conjecture has many applications in probability theory, geometric analysis and computer science.  Keller and Klein \cite{KK22} completely solved Tomaszewski's conjecture. We refer the reader to Keller and  Klein \cite{KK22} for
details.

Motivated by Tomaszewski's conjecture, the quantity
$$
\inf_{\alpha}\P\left(|X_{\alpha}-\E[X_{\alpha}]|\leq \sqrt{{\rm Var}(X_{\alpha})}\right)
$$
was studied in a series of papers (\cite{SHS23a,SHS23b,SHS23c,ZHS23}), where $\{X_{\alpha}\}$ is  a family of random variables following certain types of distributions with finite expectation $\E[X_{\alpha}]$
and finite variance $\Var(X_{\alpha})$.

Hitczenko-Kwapie\'{n}'s conjecture is a counterpart of Tomaszewski's conjecture, see \cite{HK},  and it
can be
stated as $\P(|X|\geq 1)\geq 7/32$ for all $X\in\mathcal{R}$.  Hitczenko-Kwapie\'{n}'s conjecture was recently proved by Hollom and Portier \cite{HP23}. In fact, Hollom and Portier \cite{HP23} found an explicit expression for the anti-concentration function
$$
A(y)=\inf_{X\in\mathcal{R}} \P\left(|X|\geq y\right), \quad y>0.
$$
For earlier work on Hitczenko-Kwapie\'{n}'s conjecture,
see \cite{DK22} and the references therein.

Motivated by \cite{HP23},
we study in this paper the anti-concentration function
$$
A(y)=\inf_{\alpha}\P\left(|X_\alpha-\E[X_\alpha]|\geq y\sqrt{{\rm Var}(X)}\right), \quad y>0,
$$
and find its explicit expression, where $\{X_\alpha\}$ is a family of random variables following certain familiar types of distributions.

In Section 2, we will show that, in the cases of uniform, exponential, non-degenerate Gaussian distributions and student's $t$-distribution, the anti-concentration function is not identically zero, which means that the corresponding families of random variables
have some sort of  anti-concentration property. In Section 3, we show that, in the cases of binomial, Poisson, negative binomial,
hypergeometric, Gamma, Pareto, Weibull, log-normal and Beta distributions,
the anti-concentration function is identically zero and so the corresponding families of random variables do not have the anti-concentration property above.

For information on  the distributions used in this paper, we refer the reader
to the handbook of Krishnamoorthy \cite{K}.

\section{Some distributions with anti-concentration property}\setcounter{equation}{0}

In this section, we will show that, in the cases of uniform, exponential, non-degenerate Gaussian distributions and student's $t$-distribution, the anti-concentration function is not identically zero.

\subsection{Uniform distribution on $[a,b]$ with $a<b$}

\begin{pro}\label{thmuniform}
	Let $U_{a,b}$ denote a uniform random variable on the interval $[a,b]$ with $a<b$.
	Then for $y>0$, it holds that
\begin{align}\label{pro-2.1-a}
		A_1(y):=
	\inf_{a,b}\P\left(|U_{a,b}-\E[U_{a,b}]|\ge y\sqrt{{\rm Var}(U_{a,b})}\right)
=\left\{
\begin{array}{cl}
		1-\frac{y}{\sqrt{3}},&{\rm if}\ 0<y<\sqrt{3},\\
		0,&{\rm if}\ y\ge \sqrt{3}.
\end{array}\right.
\end{align}
\end{pro}

\noindent {\bf Proof.}\ Without loss of generality,  we can assume  that $b>0,a=-b$. Then we have
\begin{align*}
A_1(y)
&=\inf_{b}\P\left(|U_{-b,b}|\ge y\sqrt{{\rm Var}(U_{-b,b})}\right)\\
&=\inf_{b}\P\left(|U_{-b,b}|\ge \frac{yb}{\sqrt{3}}\right)\\
&=\left\{
\begin{array}{cl}
		1-\frac{y}{\sqrt{3}},&{\rm if}\ 0<y<\sqrt{3},\\
		0,&{\rm if}\ y\ge \sqrt{3},
\end{array}\right.
\end{align*}
i.e. (\ref{pro-2.1-a}) holds.\hfill\fbox

\subsection{Exponential distribution}

\begin{pro}\label{thmexponential}
Let $X_{\lambda}$ denote an exponential random variable with parameter $\lambda>0$. Then it holds that
\begin{align}\label{thmexponential-1}
		A_2(y)
		:=\inf_{\lambda>0}\P\left(|X_{\lambda}-\E[X_{\lambda}]|\ge y\sqrt{{\rm Var}(X_{\lambda})}\right)=\left\{
\begin{array}{cl}
	1-e^{-(1-y)}+e^{-(1+y)},& {\rm if}\ 0<y<1,\\
			e^{-(1+y)},&{\rm if}\ y\ge 1.
\end{array}\right.
\end{align}
\end{pro}

\noindent Proof.\ \ We know that $\E[X_{\lambda}]=\frac{1}{\lambda}$ and  ${\rm Var}(X_{\lambda})=\frac{1}{\lambda^2}$. Then
$$
A_2(y)
=\inf_{\lambda>0}\P\left(\left|X_{\lambda}-\frac{1}{\lambda}\right|\ge \frac{y}{\lambda}\right).
$$
For $y\ge 1$, we have
$$
A_2(y)
=\inf_{\lambda>0}\P\left(X_{\lambda}\ge \frac{1+y}{\lambda}\right)=
\P(X_{1}\ge 1+y)=e^{-(1+y)}.$$
For $0<y<1$, we have
\begin{align*}
A_2(y)
&=\inf_{\lambda>0}\left[\P\left(X_{\lambda}\ge \frac{1+y}{\lambda}\right)+\P\left(X_{\lambda}\le \frac{1-y}{\lambda}\right)\right]\\
&=\P(X_{1}\ge 1+y)+\P(X_{1}\le1-y)\\
&=1-e^{-(1-y)}+e^{-(1+y)}.
\end{align*}
 Hence (\ref{thmexponential-1}) holds. \hfill\fbox

\subsection{Non-degenerate Gaussian distribution}

\begin{pro}\label{thmnormal}
	Let $N_{\mu,\sigma}$ be a Gaussian random variable with
	mean $\mu$ and variance $\sigma^2>0$. Then for $y>0$, it holds that
$$
A_3(y):=\P\left(|N_{\mu,\sigma}-\mu|\ge y\sigma\right)=2\Phi (-y),
$$
where $\Phi(\cdot)$ stands for the standard normal distribution function.
\end{pro}

\noindent {\bf Proof.} \  Denote $Z:=\frac{N_{\mu,\sigma}-\mu}{\sigma}$. Then $Z\sim N(0,1)$ and thus
\begin{align*}
\P\left(|N_{\mu,\sigma}-\mu|\ge y\sigma\right)&=\P\left(\left|\frac{N_{\mu,\sigma}-\mu}{\sigma}\right|\ge y\right)=P\left (|Z| \ge y \right)=2\Phi (-y).
\end{align*}
\hfill\fbox

\subsection{Student's $t$-distribution}

For $n\in \mathbb{N}_+$, let $X_n$ be a $t$-random variable with parameter $n$.
Denote by $F(a,b;c;z)$ the hypergeometric function (cf. \cite{R1}):
$$
F(a,b;c;z)=\sum_{j=0}^{\infty} \frac{(a)_j(b)_j}{(c)_j}\cdot \frac{z^j}{j!},\ \ \  \left | z \right | <1,
$$
where $(\alpha)_j\coloneqq \alpha(\alpha+1)\cdots (\alpha+j-1)$ for $j\ge 1$, and $(\alpha)_0=1$ for $\alpha\ne0$. The probability density function and the cumulative distribution function of $X_n$ are
\begin{align}
	f_n(x)&=\frac{\Gamma\left(\frac{n+1}{2}\right)}{\sqrt{n \pi} \Gamma\left(\frac{n}{2}\right)}\left(1+\frac{x^2}{n}\right)^{-\frac{n+1}{2}},\ \ x\in\mathbb{R},\label{2.4-a}\\
F_n(x)&=\frac{1}{2}+x \frac{\Gamma\left(\frac{n+1}{2}\right) }{\sqrt{n\pi}\Gamma(\frac{n}{2})}
	F\left(\frac{1}{2}, \frac{n+1}{2} ; \frac{3}{2} ;-\frac{x^2}{n}\right),\ \ x\in\mathbb{R}. \label{2.4-b}
\end{align}
The expectation $\E[X_n]=0$, and the variance  $\Var(X_n)=\frac{n}{n-2}$ for $n\ge3$.

\begin{thm}\label{thmstu}
For $y\in (0, \sqrt{6}/2)$, it holds that
\begin{align}\label{thmstu-a}
A_4(y):=\inf_{n\ge3}\P\left(|X_{n}-\E[X_{n}]|\ge y\sqrt{{\rm Var}(X_{n})}\right)=2-2\max_{3\le n\le n_0(y)+1}F_n\left(y\sqrt{\frac{n}{n-2}}\right),
\end{align}
where $n_0(y):=\min\{n\geq 3: y^2<\frac{3n^2-14n+16}{2n^2-6n+3}\}.$
\end{thm}

\noindent{\bf Proof.} \ We have
\begin{align}\label{thmstu-a}
A_4(y)
=1- \sup_{n\ge3}\P\left(|X_{n}|< y\sqrt{\frac{n}{n-2}}\right),
\end{align}
and
\begin{align}\label{thmstu-b}
\P\left(\left|X_n\right| <y \sqrt{\frac{n}{n-2}}\right) & =2y \sqrt{\frac{n}{n-2}} \frac{\Gamma\left(\frac{n+1}{2}\right)}{\sqrt{n \pi} \Gamma\left(\frac{n}{2}\right)} F\left(\frac{1}{2}, \frac{n+1}{2} ; \frac{3}{2} ;-\frac{y^2}{n-2}\right):=J_n.
\end{align}
Then we have
\begin{align}\label{thmstu-c}
\frac{J_{n+2}}{J_n}<1 &
\Leftrightarrow \frac{(n+1)(n-2)^{\frac{1}{2}}}{n^{\frac{3}{2}}} \frac{F\left(\frac{1}{2}, \frac{n+3}{2}; \frac{3}{2};-\frac{y^2}{n}\right)}{F\left(\frac{1}{2}, \frac{n+1}{2} ; \frac{3}{2};-\frac{y^2}{n-2}\right)}<1 \nonumber\\
& \Leftrightarrow F\left(\frac{1}{2}, \frac{n+3}{2} ; \frac{3}{2} ;-\frac{y^2}{n}\right)<\frac{n^{\frac{3}{2}}}{(n+1)(n-2)^{\frac{1}{2}}} F\left(\frac{1}{2}, \frac{n+1}{2}; \frac{3}{2};-\frac{y^2}{n-2}\right).
\end{align}

By the Gauss relation between contiguous functions (cf. \cite[p. 71, 21(13)]{R1})
and the fact that $F(a,b;a;z)=(1-z)^{-b}$, we get
\begin{align*}
	\frac{n+1}{2} F\left(\frac{1}{2}, \frac{n+3}{2} ; \frac{3}{2} ;-\frac{y^2}{n}\right) & =\frac{n}{2} F\left(\frac{1}{2}, \frac{n+1}{2} ; \frac{3}{2} ;-\frac{y^2}{n}\right)+\frac{1}{2} F\left(\frac{1}{2}, \frac{n+1}{2} ; \frac{1}{2} ;-\frac{y^2}{n}\right) \\ & =\frac{n}{2} F\left(\frac{1}{2}, \frac{n+1}{2} ; \frac{3}{2} ;-\frac{y^2}{n}\right)+\frac{1}{2}\left(\frac{n}{n+y^2}\right)^{\frac{n+1}{2}},
\end{align*}
which implies that
$$
F\left(\frac{1}{2}, \frac{n+3}{2} ; \frac{3}{2} ;-\frac{y^2}{n}\right)=\frac{n}{n+1} F\left(\frac{1}{2}, \frac{n+1}{2} ; \frac{3}{2} ;-\frac{y^2}{n}\right)+\frac{1}{n+1}\left(\frac{n}{n+y^2}\right)^{\frac{n+1}{2}}.
$$
Then by (\ref{thmstu-c}), (\ref{2.4-a}) and  (\ref{2.4-b}), we have
\begin{align}\label{thmstu-d}
&\frac{J_{n+2}}{J_n}<1 \nonumber\\
 \Leftrightarrow& \frac{n}{n+1} F\left(\frac{1}{2}, \frac{n+1}{2} ; \frac{3}{2} ;-\frac{y^2}{n}\right)+\frac{1}{n+1}\left(\frac{n}{n+y^2}\right)^{\frac{n+1}{2}}<\frac{n^{\frac{3}{2}}}{(n+1)(n-2)^{\frac{1}{2}}} F\left(\frac{1}{2}, \frac{n+1}{2} ; \frac{3}{2} ;-\frac{y^2}{n-2}\right) \nonumber\\
	\Leftrightarrow& F\left(\frac{1}{2}, \frac{n+1}{2} ; \frac{3}{2} ;-\frac{y^2}{n}\right)+\frac{1}{n}\left(\frac{n}{n+y^2}\right)^{\frac{n+1}{2}}<
\left(\frac{n}{n-2}\right)^{\frac{1}{2}} F\left(\frac{1}{2}, \frac{n+1}{2} ; \frac{3}{2} ;-\frac{y^2}{n-2}\right)\nonumber\\
\Leftrightarrow&\frac{\Gamma(\frac{n+1}{2})y}{\sqrt{n\pi}\Gamma(\frac{n}{2})} \left[F\left(\frac{1}{2}, \frac{n+1}{2} ; \frac{3}{2} ;-\frac{y^2}{n}\right)+\frac{1}{n}\left(\frac{n}{n+y^2}\right)^{\frac{n+1}{2}}\right]\nonumber\\
&\quad <\frac{\Gamma(\frac{n+1}{2})y}{\sqrt{n\pi}\Gamma(\frac{n}{2})}\cdot
\left(\frac{n}{n-2}\right)^{\frac{1}{2}} F\left(\frac{1}{2}, \frac{n+1}{2} ; \frac{3}{2} ;-\frac{y^2}{n-2}\right)\nonumber\\
	\Leftrightarrow& \int_{0}^{y} f_n(x)dx+\frac{\Gamma(\frac{n+1}{2})}{\sqrt{n\pi}\Gamma(\frac{n}{2})}
\cdot\frac{y}{n}\left(\frac{n}{n+y^2}\right)^{\frac{n+1}{2}}
<\int_{0}^{y(\frac{n}{n-2})^{\frac{1}{2}}}f_n(x)dx\nonumber\\
	\Leftrightarrow& \frac{y}{n}\left(\frac{n}{n+y^2}\right)^{\frac{n+1}{2}}<
\int_{y}^{y(\frac{n}{n-2})^{\frac{1}{2}}}\left(1+\frac{x^2}{n}\right)^{-\frac{n+1}{2}}dx \nonumber\\
	\Leftrightarrow&y<n\int_{y}^{y\left(\frac{n}{n-2}\right)^{\frac{1}{2}}}
\left(\frac{n+y^2}{n+x^2}\right)^{\frac{n+1}{2}}dx\nonumber\\
	\Leftrightarrow&
	1<n\int_{1}^{(\frac{n}{n-2})^{\frac{1}{2}}}
\left(\frac{n+y^2}{n+y^2t^2}\right)^{\frac{n+1}{2}}dt.
\end{align}

By the monotonicity of the function $(\frac{n+y^2}{n+y^2t^2})^{\frac{n+1}{2}}$ and Taylor's formula, we have
\begin{align*}
1<n\int_{1}^{(\frac{n}{n-2})^{\frac{1}{2}}}
\left(\frac{n+y^2}{n+y^2t^2}\right)^{\frac{n+1}{2}}dt
&\Leftarrow 1<n\left[\left(\frac{n}{n-2}\right)^{\frac{1}{2}}-1\right]
\left[\frac{(n+y^2)(n-2)}{n(n+y^2-2)}\right]^{\frac{n+1}{2}} \\
 &\Leftrightarrow
1<n\left[\left(1+\frac{2}{n-2}\right)^{\frac{1}{2}}-1\right]
\left[1-\frac{2y^2}{n(n+y^2-2)}\right]^{\frac{n+1}{2}}\\
&\Leftarrow 1<n\left[\frac{1}{n-2}-\frac{1}{2(n-2)^2}\right]
\left[1-\frac{y^2(n+1)}{n(n+y^2-2)}\right]\\
&\Leftrightarrow 1<\frac{(2n-5)(n^2-2n-y^2)}{2(n-2)^2(n+y^2-2)} \\
&\Leftrightarrow (2y^2-3)n^2+(14-6y^2)n+3y^2-16<0.
\end{align*}
When $0<y<\frac{\sqrt{6}}{2}$, we have  $y^2<\frac{3}{2}$, and thus
$$
(2y^2-3)n^2+(14-6y^2)n+3y^2-16<0 \Leftrightarrow y^2<\frac{3n^2-14n+16}{2n^2-6n+3}.
$$
Note that
$
\lim_{n\to\infty}\frac{3n^2-14n+16}{2n^2-6n+3}=\frac{3}{2},
$
and  $\{\frac{3n^2-14n+16}{2n^2-6n+3},n\geq 3\}$ is an increasing sequence. Thus
we know that for any $0<y<\frac{\sqrt{6}}{2}$, the set $\{n\geq 3: y^2<\frac{3n^2-14n+16}{2n^2-6n+3}\}$ is not empty. Define
$$
n_0(y):=\min\left\{n\geq 3: y^2<\frac{3n^2-14n+16}{2n^2-6n+3}\right\}.$$
Then by the symmetry of the density function $f_n(x)$, we have
\begin{align*}
A_4(y)
&=1-\sup_{n\geq 3}J_n=1-\max_{3\leq n\leq n_0(y)+1}J_n=
1-\max_{3\leq n\leq n_0(y)+1}\P\left(\left|X_n\right| <y \sqrt{\frac{n}{n-2}}\right)\\
&=1-\max_{3\leq n\leq n_0(y)+1}\left[F_n\left(y\sqrt{\frac{n}{n-2}}\right)
-F_n\left(-y\sqrt{\frac{n}{n-2}}\right)\right]\\
&=1-\max_{3\leq n\leq n_0(y)+1}\left[F_n\left(y\sqrt{\frac{n}{n-2}}\right)
-\left(1-F_n\left(y\sqrt{\frac{n}{n-2}}\right)\right)\right]\\
&=2-2\max_{3\leq n\leq n_0(y)+1}F_n\left(y\sqrt{\frac{n}{n-2}}\right).
\end{align*}
The proof is complete.\hfill\fbox

\begin{rem}
If $0<y\leq 1$ and $t\geq 1$, by $\frac{n+y^2}{n+y^2t^2}=\frac{\frac{n}{y^2}+1}{\frac{n}{y^2}+t^2}
=1-\frac{t^2-1}{\frac{n}{y^2}+t^2}\geq 1-\frac{t^2-1}{n+t^2}=\frac{n+1}{n+t^2}$,  we know that
	\begin{align*}
		n\int_{1}^{(\frac{n}{n-2})^{\frac{1}{2}}}
		\left(\frac{n+y^2}{n+y^2t^2}\right)^{\frac{n+1}{2}}dt&\geq n\int_{1}^{(\frac{n}{n-2})^{\frac{1}{2}}}
		\left(\frac{n+1}{n+t^2}\right)^{\frac{n+1}{2}}dt.
	\end{align*}
	By the proof of \cite[Theorem 3.1]{SHS23b}, we know that
	$$
	n\int_{1}^{(\frac{n}{n-2})^{\frac{1}{2}}}
	\left(\frac{n+1}{n+t^2}\right)^{\frac{n+1}{2}}dt>1,\ \ \ \forall n\geq3.
	$$
	It follows that, for $y\in (0, 1]$, we have
	$$
	n\int_{1}^{(\frac{n}{n-2})^{\frac{1}{2}}}
	\left(\frac{n+y^2}{n+y^2t^2}\right)^{\frac{n+1}{2}}dt>1,
	$$
	which together with (\ref{thmstu-d}) implies that, when $y\in (0, 1]$, for any $n\geq 3$, it holds that $J_{n+2}<J_n$.
	Hence for $y\in (0, 1]$,
	\begin{align*}
		A_4(y)
		&=1-\max_{n=3,4}J_n=2-2\max_{n=3,4}F_n\left(y\sqrt{\frac{n}{n-2}}\right).
	\end{align*}

\end{rem}

\section{Some distributions without anti-concentration property}\setcounter{equation}{0}

In this section, we will show that, in the cases of binomial, Poisson, negative binomial, hypergeometric, Gamma, Pareto, Weibull, log-normal and Beta distributions, the anti-concentration function is identically zero.

\subsection{Binomial distribution}

\begin{pro}\label{thmbin}
	Let $B(n, p)$ denote a binomial random variable with parameters $n$ and $p$. For any $y>0$, we have
\begin{align}\label{thmbin-a}
A_5(y)
:=\inf_{0<p\leq 1}\P\left(|B_{n,p}-\E[B_{n,p}]|\ge y\sqrt{{\rm Var}(B_{n,p})}\right)=0.
\end{align}
\end{pro}	

\noindent{\bf  Proof.}\ We have
$A_5(y)=\inf_{0<p\leq 1}\P(|B_{n,p}-np|\ge y\sqrt{np(1-p)})$.  It is  easy to see
that $np-y\sqrt{np(1-p)}<0$ if and only if $p<\frac{y^2}{n+y^2}$. In addition,   $\lim\limits_{p\downarrow 0}(np+y\sqrt{np(1-p)})=0$.
Thus for $y>0$, we have
$$
\lim_{p\downarrow 0}\P\left(|B_{n,p}-np|\ge y\sqrt{np(1-p)}\right)=\lim_{p\downarrow 0}(1-\P(B_{n,p}=0))=\lim_{p\downarrow 0}(1-(1-p)^n)= 0,
$$
hence (\ref{thmbin-a}) holds.\hfill\fbox

\subsection{Poisson distribution}

 \begin{pro}\label{thmpoi2}
  Let $N_\lambda$  denote a Poisson random variable with parameter $\lambda>0$.  For any $y>0$,  we have
 \begin{align}
  A_6(y):=\inf_{\lambda>0}\P\left(|N_\lambda-\E[N_\lambda]|\ge y\sqrt{{\rm Var}(N_\lambda)}\right)=0.
 \end{align}
 \end{pro}

\noindent{\bf Proof.} Recall that $\E[N_\lambda]={\rm Var}(N_\lambda)=\lambda$.
When $0<\lambda<y^2$, we have $\lambda-y\sqrt{\lambda}<0 $, and thus
$$
\P\left(|N_\lambda-\E[N_\lambda]|\ge y\sqrt{{\rm Var}(N_\lambda)}\right)=\P(N_\lambda\ge \lambda+y\sqrt{\lambda}).
$$
Note that $0<\sqrt{\lambda}<\frac{-y+\sqrt{y^2+4}}{2}$ implies that $0<\lambda+y\sqrt{\lambda}\le 1$. So when $0<\lambda<\min\{y^2,(\frac{-y+\sqrt{y^2+4}}{2})^2\}$, we have
$$
\P\left(N_\lambda\ge \lambda+y\sqrt{\lambda}\right)=1-\P(N_\lambda=0)=1-e^{-\lambda}.
$$
It follows that $A_6(y)=0$.\hfill\fbox

\subsection{Negative binomial distribution}

\begin{pro}\label{thmneg}
For $r>0$ and $p\in (0, 1]$, let $X_{r,p}$ denote a negative binomial random variable with parameters $(r, p)$, that is,
 \begin{eqnarray}\label{1.1-b}
P(X_{r,p}=l)=\binom{-r}{l}p^r(-q)^l=\binom{r+l-1}{l}p^rq^{l},\ l=0,1,2,\ldots,
\end{eqnarray}
where $q=1-p$.
 For any $y>0$, we have
\begin{align}\label{thmneg-a}
A_7(y)
:=\inf_{0<p\leq 1}\P\left(|X_{r,p}-\E[X_{r,p}]|\ge y\sqrt{{\rm Var}(X_{r,p})}\right)=0.
\end{align}
\end{pro}

\noindent {\bf Proof.}\
Recall that $\E[X_{r,p}]=\frac{rq}{p}$ and the variance is $\Var(X_{r,p})=\frac{rq}{p^2}$.
Thus
$$
A_7(y)=\inf_{0<p\leq 1}\P\left(|X_{r,p}-\frac{rq}{p}|\ge y\sqrt{\frac{rq}{p^2}}\right).
$$
We consider the situation when $\frac{rq}{p}-y\sqrt{\frac{rq}{p^2}}< 0$, i.e., $p\in(1-\frac{y^2}{r},1]$. We know  $\lim\limits_{p\uparrow 1}\frac{rq}{p}+y\sqrt{\frac{rq}{p^2}}=0$. It follows that
$$
\lim_{p\uparrow 1}\P\left(|X_{p}-\E[X_{p}]|\ge y\sqrt{{\rm Var}(X_{p})}\right)=\lim_{p\uparrow 1}\left(1-\P(X_p=0)\right)=\lim_{p\uparrow 1}\left(1-p^r\right)=0.
$$
Hence (\ref{thmneg-a}) holds.\hfill\fbox

\subsection{Hypergeometric distribution}

\begin{pro}\label{thmhyper}
For $M,N,n \in\mathbb{Z_+}$ satisfying that $M\leq N, n\leq N$, let $X_{M,N,n}$ be a hypergeometric random variable with parameters $M,N,n \in\mathbb{Z_+}$, that is,
$$
\P(X_{M,N,n}=k)=\frac{\binom{M}{k}\binom{N-M}{n-k}}{\binom{N}{n}}, \quad k=0,1,\cdots,\min(M,n).
$$
 For any $y>0$, we have
\begin{align}\label{thmhyper-a}
A_8(y):=\inf\P\left(|X_{M,N,n}-E[X_{M,N,n}]|\ge y\sqrt{{\rm Var}(X_{M,N,n})}\right)=0,
\end{align}
where then infimum is taken over all $M,N,n \in\mathbb{Z_+}$ satisfying that $M\leq N, n\leq N$.
\end{pro}

\noindent {\bf Proof.} \
Recall that $\E[X_{M,N,n}]=\frac{nM}{N}$ and $\Var(X_{M,N,n})=\frac{nM}{N}(1-\frac{M}{N})\cdot \frac{N-n}{N-1}$.
We consider the situation when $M=N-1$ and $n=1$. Then
$$
\P(X_{N-1,N,1}=k)=\left\{
\begin{array}{cl}
	\frac{1}{N},& {\rm if}\ k=0,\\
	\frac{N-1}{N},&{\rm if}\ k=1.
\end{array}\right.
$$
We have  $\E[X_{N-1,N,1}]=\frac{N-1}{N}$ and  $\Var(X_{N-1,N,1})=\frac{N-1}{N^2}$, and thus
\begin{align*}
&\P\left(|X_{N-1,N,1}-\E[X_{N-1,N,1}]|\ge y\sqrt{{\rm Var}(X_{N-1,N,1})}\right)\\
&=\P\left(X_{N-1,N,1}\le\frac{N-1}{N}-y\sqrt{\frac{N-1}{N^2}}\right)
+\P\left(X_{N-1,N,1}\ge\frac{N-1}{N}+y\sqrt{\frac{N-1}{N^2}}\right).
\end{align*}
When $N$ is large enough, we have $0<\frac{N-1}{N}-y\sqrt{\frac{N-1}{N^2}}<1$ and $\frac{N-1}{N}+y\sqrt{\frac{N-1}{N^2}}>1$, and thus
$$
\P\left(|X_{N-1,N,1}-E[X_{N-1,N,1}]|\ge y\sqrt{{\rm Var}(X_{N-1,N,1})}\right)=\P\left(X_{N-1,N,1}\le\frac{N-1}{N}-y\sqrt{\frac{N-1}{N^2}}\right)
=\frac{1}{N}.
$$
Hence (\ref{thmhyper-a}) holds.\hfill\fbox

\subsection{Gamma distribution}

\begin{pro}\label{thmgamma}	
For $\alpha>0,\beta>0$, let $X_{\alpha,\beta}$ denote a Gamma random variable with probability density function
$$
f_{\alpha,\beta}(x)=\frac{x^{\alpha-1}e^{-x/\beta}}{\Gamma(\alpha)\beta^\alpha},\ \ x>0.
$$
For any $y>0$, we have
\begin{align}
A_9(y):=\inf_{\alpha>0,\beta>0}\P\left(|X_{\alpha,\beta}-\E[X_{\alpha,\beta}]|\ge y\sqrt{{\rm Var}(X_{\alpha,\beta})}\right)=0.
\end{align}
\end{pro}

\noindent {\bf Proof.} \ Recall that $\E[X_{\alpha,\beta}]=\alpha\beta$ and $\Var(X_{\alpha,\beta})=\alpha \beta^2$.
Thus
$$
A_9(y)=\inf_{\alpha>0,\beta>0}\P\left(|X_{\alpha,\beta}-\alpha\beta|\ge y\sqrt{\alpha}\beta\right).
$$
Set $g(\alpha):=\P\left\{|X_{\alpha,1}-\alpha|\ge y\sqrt{\alpha}\right\}$. In the following, we will prove that for $\forall y>0$, $\lim_{\alpha \to 0}g(\alpha)=0 $, thus we get $A_9(y)=0$.

\noindent For any $y>0$,
  if $0<\alpha<y^2$, then $\alpha-y\sqrt{\alpha}<0$.
 Then $g(\alpha)=\P(X_{\alpha,1}\ge\alpha+y\sqrt{\alpha})$ and thus
\begin{align*}
g(\alpha)&=1-\int_{0}^{\alpha+y\sqrt{\alpha}} \frac{x^{\alpha-1}e^{-x}}{\Gamma(\alpha)}dx\le 1-\int_{0}^{\alpha} \frac{x^{\alpha-1}e^{-x}}{\Gamma(\alpha)}dx\\
&\le1-\frac{e^{-\alpha}}{\Gamma(\alpha)}\int_{0}^{\alpha}x^{\alpha-1}dx
=1-e^{-\alpha}\frac{\alpha^{\alpha-1}}{\Gamma(\alpha)}.
\end{align*}

\noindent Then we have $$
\limsup\limits_{\alpha \to 0} g(\alpha)\le1-\lim_{\alpha \to 0}\frac{\alpha^{\alpha-1}}{\Gamma(\alpha)}=1-\lim_{\alpha \to 0}\frac{\alpha^\alpha}{\Gamma(\alpha+1)}=1-\lim_{\alpha \to 0}e^{\alpha\ln \alpha}=0.
 $$
Thus we have $\lim_{\alpha \to 0}g(\alpha)=0$. The proof is complete.\hfill\fbox

\subsection{Pareto distribution}

\begin{pro}\label{thmpareto}
For $r>0,A>0$, let $X_{r,A}$ denote a Pareto random variable with probability density function
$$
f_{r,A}=rA^r\frac{1}{x^{r+1}},\ \ x\ge A.
$$
For any $y>0$, we have
\begin{align}\label{thmpareto-a}
A_{10}(y):=\inf_{r>2,A>0}\P\left(|X_{r,A}-\E[X_{r,A}]|\ge y\sqrt{{\rm Var}(X_{r,A})}\right)=0.
\end{align}
\end{pro}

\noindent {\bf Proof.} \ Recall that $\E[X_{r,A}]=\frac{r}{r-1}A$ when $r>1$ and $\Var(X_{r,A})=\frac{rA^2}{(r-2)(r-1)^2}$ when $r>2$.
Thus $A_{10}(y)=\inf_{r>2,A>0}\P(|X_{r,A}-\frac{r}{r-1}A|\ge \frac{yA}{r-1}\sqrt{\frac{r}{r-2}})$. For any $y>0$, if $r$ is close to 2 from above, then $\sqrt{\frac{r-2}{r}}\le y$ and thus $\frac{r}{r-1}A-\frac{yA}{r-1}\sqrt{\frac{r}{r-2}}\le A$, and so $\P(|X_{r,A}-\frac{r}{r-1}A|\ge \frac{yA}{r-1}\sqrt{\frac{r}{r-2}})
=\P( X_{r,A}\ge \frac{r}{r-1}A+\frac{yA}{r-1}\sqrt{\frac{r}{r-2}})$. Hence  it is enough to show that
\begin{align}\label{thmpareto-b}
\lim_{r\downarrow 2} \P\left(X_{r,A}\ge \frac{r}{r-1}A+\frac{yA}{r-1}\sqrt{\frac{r}{r-2}}\right)=0.
\end{align}
Obviously, we have
$\frac{r}{r-1}A+\frac{yA}{r-1}\sqrt{\frac{r}{r-2}}>A$.  Thus
\begin{align*}
&\lim_{r\downarrow 2} \P\left(X_{r,A}\ge \frac{r}{r-1}A+\frac{yA}{r-1}\sqrt{\frac{r}{r-2}}\right)\\
&=\lim_{r\downarrow 2} \int_{\frac{r}{r-1}A+\frac{yA}{r-1}\sqrt{\frac{r}{r-2}}}^{\infty}rA^r\frac{1}{x^{r+1}} dx\\
&=\lim_{r\to  2}\int_{\frac{r}{r-1}+\frac{y}{r-1}\sqrt{\frac{r}{r-2}}}^{\infty}r\frac{1}{y^{r+1}} dy  \\
&=\lim_{r\to 2}\left(\frac{r-1}{r+y\sqrt{\frac{r}{r-2}}}\right)^r=0.
\end{align*}
Hence (\ref{thmpareto-b}) holds and the proof is complete.\hfill\fbox

\subsection{Weibull distribution}

\begin{pro}\label{thmweibull}	
For $\alpha>0,\lambda>0$, let $X_{\alpha,\lambda}$ denote a Weibull random variable with probability density function
$$
f_{\alpha,\lambda}=\alpha\lambda x^{\alpha-1}e^{-\lambda x^\alpha},\ \ x>0.
$$
For any $y>0$, we have
\begin{align}\label{thmweibull-a}
A_{11}(y)
:=\inf_{\alpha>0,\lambda>0}\P\left(|X_{\alpha,\lambda}-\E[X_{\alpha,\lambda}]|\ge y\sqrt{{\rm Var}(X_{\alpha,\lambda})}\right)=0.
\end{align}
\end{pro}

\noindent {\bf Proof.} \
Recall that $\E[X_{\alpha,\lambda}]=\lambda^{-\frac{1}{\alpha}}\Gamma(1+\frac{1}{\alpha})$ and
$\Var(X_{\alpha,\lambda})=\lambda^{-\frac{2}{\alpha}}
[\Gamma(1+\frac{2}{\alpha})-(\Gamma(1+\frac{1}{\alpha}))^2]$. Thus
$$F_7(y)=\inf_{\alpha>0,\lambda>0}\P\left(\left|X_{\alpha,\lambda}-
\lambda^{-\frac{1}{\alpha}}\Gamma(1+\frac{1}{\alpha})\right|\ge y\lambda^{-\frac{1}{\alpha}}\sqrt{\Gamma(1+\frac{2}{\alpha})-(\Gamma(1+\frac{1}{\alpha}))^2}
\right).
$$
Firstly, by the Stirling formula and the monotonicity of the Gamma function, we prove that $\lim_{\alpha\to 0}\frac{\Gamma(1+\frac{2}{\alpha})}{(\Gamma(1+\frac{1}{\alpha}))^2}=+\infty$:
\begin{align*}
	\liminf\limits_{\alpha\to 0}\frac{\Gamma(1+\frac{2}{\alpha})}{(\Gamma(1+\frac{1}{\alpha}))^2}&=	\liminf\limits_{x\to +\infty}\frac{\Gamma(1+2x)}{(\Gamma(1+x))^2}\ge \liminf\limits_{x\to +\infty} \frac{\Gamma(1+2[x])}{(\Gamma(2+[x]))^2}\\
	&=\liminf\limits_{n\to +\infty} \frac{\Gamma(2n+1)}{(\Gamma(n+2))^2}=\lim_{n\to +\infty}\frac{(2n)!}{(n+1)!(n+1)!}\\&=\lim_{n\to +\infty}\frac{\sqrt{2\pi}e^{-2n}(2n)^{2n+\frac{1}{2}}}{2\pi e^{-(2n+2)}(n+1)^{2n+3}}\\&=\frac{e^2}{\sqrt{2\pi}}\lim_{n\to +\infty}
	\frac{2^{2n+\frac{1}{2}}}{(1+\frac{1}{n})^{2n+\frac{1}{2}}(n+1)^{\frac{5}{2}}}=+\infty.
\end{align*}
Therefore, for any $y>0$, as $\alpha\to 0$, we have
\begin{align*}
(1+y^2)(\Gamma(1+\frac{1}{\alpha}))^2-y^2\Gamma(1+\frac{2}{\alpha})
=(\Gamma(1+\frac{1}{\alpha}))^2\left[
(1+y^2)-y^2\frac{\Gamma(1+\frac{2}{\alpha})}{(\Gamma(1+\frac{1}{\alpha}))^2}\right]\to -\infty,
\end{align*}
which implies that as $\alpha>0$ and is small enough,
\begin{align*}
\Gamma(1+\frac{1}{\alpha})-y\sqrt{\Gamma(1+\frac{2}{\alpha})-
(\Gamma(1+\frac{1}{\alpha}))^2}=\frac{(1+y^2)(\Gamma(1+\frac{1}{\alpha}))^2-y^2\Gamma(1+\frac{2}{\alpha})}{\Gamma(1+\frac{1}{\alpha})+y\sqrt{\Gamma(1+\frac{2}{\alpha})-
(\Gamma(1+\frac{1}{\alpha}))^2}}<0.
\end{align*}
So when $\alpha>0$ and is small enough, we have
\begin{align*}
&\P\left(\left|X_{\alpha,1}-\Gamma(1+\frac{1}{\alpha})\right|\ge y\sqrt{\Gamma(1+\frac{2}{\alpha})-(\Gamma(1+\frac{1}{\alpha}))^2}\right)\\
&=\P\left( X_{\alpha,1}\ge \Gamma(1+\frac{1}{\alpha})+y\sqrt{\Gamma(1+\frac{2}{\alpha})-(\Gamma(1+\frac{1}{\alpha}))^2}    \right).
\end{align*}

 Define $g(\alpha):=P(X_{\alpha,1}\ge \Gamma(1+\frac{1}{\alpha})+y\sqrt{\Gamma(1+\frac{2}{\alpha})-(\Gamma(1+\frac{1}{\alpha}))^2}  )$.  Then it is enough to show that $\lim_{\alpha\to0}g(\alpha)=0.$ We have
$$
g(\alpha)\le\int_{\Gamma(1+\frac{1}{\alpha})}^{\infty} \alpha x^{\alpha-1}e^{- x^\alpha}dx=e^{-(\Gamma(1+\frac{1}{\alpha}))^\alpha},
$$
and
$$
\lim_{\alpha\to0}(\Gamma(1+\frac{1}{\alpha}))^\alpha=\lim_{\alpha\to0}e^{\alpha\ln \Gamma(1+\frac{1}{\alpha}) }=e^{\lim_{x\to\infty}\frac{\ln\Gamma(1+x)}{x}}=
e^{\lim_{x\to\infty}\frac{\Gamma'(1+x)}{\Gamma(1+x)}}.
$$
Denote $\psi(z):=\frac{\Gamma'(z)}{\Gamma(z)}$. Then we have (cf. \cite{E1})
\begin{align*}
\psi(1+n)&=1+\frac{1}{2}+\frac{1}{3}+\cdots+\frac{1}{n}-\gamma,\\
\psi(z+n)&=\frac{1}{z}+\frac{1}{z+1}+\frac{1}{z+2}+\cdots+\frac{1}{z+n-1}+\psi(z)\ \ \ n=1,2,3,\cdots,
\end{align*}
where $\gamma$ is the Euler's constant. It follows that  $\lim_{x\to\infty}\frac{\Gamma'(1+x)}{\Gamma(1+x)}=+\infty$, and thus $\lim_{\alpha\to0}(\Gamma(1+\frac{1}{\alpha}))^\alpha=+\infty$.  Hence $\lim_{\alpha\to0}g(\alpha)=0$ and thus (\ref{thmweibull-a}) holds.  The proof is complete.\hfill\fbox

\subsection{Log-normal distribution}

\begin{pro}\label{thmlog}
For $\alpha,\sigma>0$, let $X_{\alpha,\sigma}$ denote a log-normal random variable with probability density function
$$
f_{\alpha,\sigma}=\frac{1}{\sigma x\sqrt{2\pi}}e^{-\frac{(\ln x-\alpha)^2}{2\sigma^2}},\ \ x>0.
$$
For any $y>0$, we have
\begin{align}\label{thmlog-a}
A_{12}(y)
:=\inf_{\alpha,\sigma>0}\P\left(|X_{\alpha,\sigma}-\E[X_{\alpha,\sigma}]|\ge y\sqrt{{\rm Var}(X_{\alpha,\sigma})}\right)=0.
\end{align}
\end{pro}

\noindent {\bf Proof.} \
Recall that $\E[X_{\alpha,\sigma}]=e^{\alpha+\sigma^2/2}$and  $\Var(X_{\alpha,\sigma})=e^{2\alpha+\sigma^2}(e^{\sigma^2}-1)$.
Define $Y:=(\ln X_{\alpha,\sigma}-a)/\sigma$. Then $Y$ is a standard normal  random variable.
Recall  $\Phi$ to denote the standard normal distribution function.

We consider the situation when $e^{\alpha+\sigma^2/2}-y\sqrt{e^{2\alpha+\sigma^2}(e^{\sigma^2}-1)}\le 0$, i.e. $\sigma^2>\ln \frac{1+y^2}{y^2}$. We have
\begin{align*}
&\inf_{\alpha,\sigma>\sqrt{\ln \frac{1+y^2}{y^2}}}\P\left(|X_{\alpha,\sigma}-e^{\alpha+\sigma^2/2}|\ge y\sqrt{e^{2\alpha+\sigma^2}(e^{\sigma^2}-1)}\right)\\
&=\inf_{\sigma>\sqrt{\ln \frac{1+y^2}{y^2}}}\P\left(Y\ge \frac{\sigma}{2}+\frac{\ln (1+y\sqrt{e^{\sigma^2}-1})}{\sigma}\right)\\
&=1-\sup_{\sigma>\sqrt{\ln \frac{1+y^2}{y^2}}}\Phi\left(\frac{\sigma}{2}+\frac{\ln (1+y\sqrt{e^{\sigma^2}-1})}{\sigma}\right)\\
&\to 0\ {\rm as}\ \sigma\to+\infty.
\end{align*}
Hence (\ref{thmlog-a}) holds.\hfill\fbox

\subsection{Beta distribution}

 \begin{pro}\label{thmbeta}
 For $p>0,q>0$, let $X_{p,q}$ denote a beta random variable with probability density function
 $$
 f_{p,q}=\frac{\Gamma(p+q)}{\Gamma(p)\Gamma(q)}x^{p-1}(1-x)^{q-1},\ \ 0<x<1.
 $$
 For any $y>0$, we have
 \begin{align}\label{thmbeta-a}
  A_{13}(y)
 :=\inf_{p>0,q>0}\P\left(|X_{p,q}-\E[X_{p,q}]|\ge y\sqrt{{\rm Var}(X_{p,q})}\right)=0.
 \end{align}
 \end{pro}

 \noindent {\bf Proof.} \
 Recall that $\E[X_{p,q}]=\frac{p}{p+q}$ and $\Var(X_{p,q})=\frac{pq}{(p+q)^2(p+q+1)}$.
 We consider the case that $p=1$.  Then we have
 \begin{align*}
 &\P\left(|X_{p,q}-\E[X_{p,q}]|\ge y\sqrt{{\rm Var}(X_{p,q})}\right)\\
 &=\P\left(\left|X_{1,q}-\frac{1}{1+q}\right|\ge y\sqrt{\frac{q}{(1+q)^2(2+q)}}\right)\\
 &=\P\left(X_{1,q}\le \frac{1}{1+q}-y\sqrt{\frac{q}{(1+q)^2(2+q)}}\right)\\
 &\quad +\P\left(X_{1,q}\ge \frac{1}{1+q}+y\sqrt{\frac{q}{(1+q)^2(2+q)}}\right).
 \end{align*}
When $q>0$ and is small enough, we have that  $0<\frac{1}{1+q}-y\sqrt{\frac{q}{(1+q)^2(2+q)}}<1$ and $\frac{1}{1+q}+y\sqrt{\frac{q}{(1+q)^2(2+q)}}>1$, and thus
 $$
  \P\left(|X_{p,q}-E[X_{p,q}]|\ge y\sqrt{{\rm Var}(X_{p,q})}\right)=\P\left(X_{1,q}\le \frac{1}{1+q}-y\sqrt{\frac{q}{(1+q)^2(2+q)}}\right).
 $$
Then it is enough to show that
$$
\lim_{q \downarrow  0}\P\left(X_{1,q}\le \frac{1}{1+q}-y\sqrt{\frac{q}{(1+q)^2(2+q)}}\right)=0.
$$
By the fact that $f_{1,q}=q(1-x)^{q-1}$, we have
\begin{align*}
&	\limsup\limits_{q \downarrow  0} \P\left(X_{1,q}\le \frac{1}{1+q}-y\sqrt{\frac{q}{(1+q)^2(2+q)}}\right)\\
& \leq
	\lim_{q \downarrow  0}\P\left(X_{1,q}\le \frac{1}{1+q}\right) \\
&=\lim_{q \downarrow  0}\int_0^{\frac{1}{1+q}}q(1-x)^{q-1}dx\\
	&=\lim_{q \downarrow  0}\left[1-\left(\frac{q}{1+q}\right)^q\right]=0.
\end{align*}
Hence (\ref{thmbeta-a}) holds and the proof is complete.\hfill\fbox

\bigskip

\noindent {\bf\large Acknowledgments}\quad  This work was supported by the National Natural Science Foundation of China
(12171335, 11931004, 12071011), the Science Development Project of Sichuan University (2020SCUNL201) and the Simons Foundation (\#960480).

\end{document}